\documentclass[thmsa]{article}
\usepackage{amssymb}
\begin{document}

\author{David Carf\`{i}}
\title{Multiplicative operators in the spaces of Schwartz families}
\date{}
\maketitle

\begin{abstract}
In this paper we introduce and study the multiplication among smooth functions and Schwartz families. This multiplication is fundamental in the formulation and development of a spectral theory for Schwartz linear operators in distribution spaces, to introduce efficiently the
Schwartz eigenfamilies of such operators and to build up a functional calculus for them. The definition of eigenfamily is absolutely natural and this new operation allows us to develop a rigorous and manageable spectral theory for Quantum Mechanics, since it appears in a form extremely similar to the current use in Physics.
\end{abstract}

\bigskip

\section{\textbf{Introduction}}

\bigskip

In the Spectral Theory of $^{\mathcal{S}}$linear operators, the eigenvalues
corresponding to the elements of certain $^{\mathcal{S}}$families have
fundamental importance. If $L$ is an $^{\mathcal{S}}$linear operator and $v$
is an $^{\mathcal{S}}$family, the family $v$ is defined an eigenfamily of
the operator $L$ if there exists a real or complex function $l$ - defined on
the set of indices of the family $v$ - such that the relation 
\[
L(v_{p})=l(p)v_{p}, 
\]
holds for every index $p$ of the family $v$. As we already have seen, in the
context of $^{\mathcal{S}}$linear operators, it is important how the
operator $L$ acts on the entire family $v$. Taking into account the above
definition, it is natural to consider the image family $L(v)$ as the product
- in pointwise sense - of the family $v$ by the function $l$, but:

\begin{itemize}
\item  \emph{is the pointwise multiplication an operation in the space of }$%
^{\mathcal{S}}$\emph{families?}

\item  \emph{what kind of properties are satisfied by this product?}
\end{itemize}

In this chapter we define and study the properties of such product.

\bigskip

\section{$^{\mathcal{O}_{M}}$\textbf{Functions}}

\bigskip

We recall, for convenience of the reader, some basic notions from theory of
distributions.

\bigskip

\textbf{Definition (of slowly increasing smooth function).}\emph{\ We denote
by }$\mathcal{O}_{M}(\Bbb{R}^{n},\Bbb{K})$\emph{, or more simply by }$%
\mathcal{O}_{M}^{(n)}$\emph{, the subspace of all smooth functions }$f$\emph{%
,\ belonging to the space }$\mathcal{C}^{\infty }(\Bbb{R}^{n},\Bbb{K})$\emph{%
, such that, for every test function }$\phi \in \mathcal{S}_{n}$\emph{\ the
product }$\phi f$\emph{\ lives in }$\mathcal{S}_{n}.$\emph{\ The space }$%
\mathcal{O}_{M}(\Bbb{R}^{n},\Bbb{K)}$\emph{\ is said to be the \textbf{space
of smooth} \textbf{functions from} }$\Bbb{R}^{n}$\emph{\ \textbf{into the
field}\ }$\Bbb{K}$\emph{\ \textbf{slowly increasing at infinity (with all
their derivatives)}.}

\bigskip

In other terms, the functions $f$ belonging to the space $\mathcal{O}%
_{M}^{(n)}$ are the only smooth functions which can generate a
multiplication operator 
\[
M_{f}:\mathcal{S}_{n}\rightarrow \mathcal{S}_{n} 
\]
of the space $\mathcal{S}_{n}$ into the space $\mathcal{S}_{n}$ itself,
(obviously) by the relation 
\[
M_{f}(g)=fg. 
\]
This is the motivation of the importance of these functions in Distribution
Theory, and the symbol itself $\mathcal{O}_{M}$ depends on this fact ($%
\mathcal{O}_{M}$ stands for multiplicative operators).

\bigskip

Let us see a first characterization.

\bigskip

\textbf{Proposition. }\emph{Let }$f\in \mathcal{E}_{n}$\emph{\ be a smooth
function. Then the following conditions are equivalent:}

\begin{itemize}
\item[\emph{1)}]  \emph{\ for all multi-index }$p\in \Bbb{N}_{0}^{n}$\emph{\
there is a polynomial }$P_{p}$\emph{\ such that, for any point }$x\in \Bbb{R}%
^{n}$\emph{, the following inequality holds} 
\[
\left| \partial ^{p}f(x)\right| \leq \left| P_{p}(x)\right| ; 
\]

\item[\emph{2)}]  \emph{\ for any test function }$\phi \in \mathcal{S}_{n}$%
\emph{\ the product }$\phi f$\emph{\ lies in }$\mathcal{S}_{n}$\emph{;}

\item[\emph{3)}]  \emph{\ for every multi-index }$p\in \Bbb{N}_{0}^{n}$\emph{%
\ and for every test function }$\phi \in \mathcal{S}_{n}$\emph{\ the product 
}$\left( \partial ^{p}f\right) \phi $\emph{\ is bounded in }$\Bbb{R}^{n}$%
\emph{.}
\end{itemize}

\bigskip

\subsection{\textbf{Topology}}

\bigskip

The standard topology of the space $\mathcal{O}_{M}^{(n\Bbb{)}}$\ is the
locally convex topology defined by the family of seminorms 
\[
\gamma _{\phi ,p}(\phi )=\sup_{x\in \Bbb{R}^{n}}\left| \phi (x)\partial
^{p}f(x)\right| 
\]
with $\phi \in \mathcal{S}_{n}$\ and $p\in \Bbb{N}_{0}^{n}$. This topology
does not have a countable basis. Also, it can be shown that the space $%
\mathcal{O}_{M}^{(n\Bbb{)}}$ is a complete space. A sequence $(f_{j})_{j\in 
\Bbb{N}}$ converges to zero in $\mathcal{O}_{M}^{(n\Bbb{)}}$ if and only if
for every test function $\phi \in \mathcal{S}_{n}$ and for every multi-index 
$p\in \Bbb{N}_{0}^{n}$, the sequence of functions $(\phi \partial
^{p}f_{j})_{j\in \Bbb{N}}$ converges to zero uniformly on $\Bbb{R}^{n}$; or,
equivalently, if, for every test function $\phi \in \mathcal{S}_{n}$, the
sequence $(\phi f_{j})_{j\in \Bbb{N}}$ converges to zero in $\mathcal{S}_{n}$%
. A filter $F$ on $\mathcal{O}_{M}^{(n\Bbb{)}}$ converges to zero in $%
\mathcal{O}_{M}^{(n\Bbb{)}}$ if and only if for every test function $\phi
\in \mathcal{S}_{n}$, the filter $\phi F$ converges to zero in $\mathcal{S}%
_{n}$.

\bigskip

\subsection{\textbf{Bounded sets in }$\mathcal{O}_{M}^{(n\Bbb{)}}$}

\bigskip

A subset $B$ of $\mathcal{O}_{M}^{(n\Bbb{)}}$ is bounded (in the topological
vector space $\mathcal{O}_{M}^{(n\Bbb{)}}$) if and only if, for all
multi-index $p\in \Bbb{N}_{0}^{n}$, there is a polynomial $P_{p}$ such that,
for any function $f\in B$, the following inequality holds true 
\[
\left| \partial ^{p}f(x)\right| \leq P_{p}(x), 
\]
for any point $x\in \Bbb{R}^{n}$.

\bigskip

\subsection{\textbf{Multiplication in }$\mathcal{S}_{n}$\textbf{\ by }$%
\mathcal{O}_{M}^{(n\Bbb{)}}$ \textbf{functions}}

\bigskip

The bilinear map 
\[
\Phi :\mathcal{O}_{M}^{(n\Bbb{)}}\times \mathcal{S}_{n}\rightarrow \mathcal{S%
}_{n}:(\phi ,f)\mapsto \phi f 
\]
is separately continuous with respect to the usual topologies of the spaces $%
\mathcal{O}_{M}^{(n\Bbb{)}}$ and $\mathcal{S}_{n}$. It follows immediately
that the multiplication operator $M_{f}$, associated with an $^{\mathcal{O}%
_{M}}$function $f$, is continuous (with respect to the standard topology of
the Schwartz space $\mathcal{S}_{n}$). Moreover, the transpose of the
operator $M_{f}$ is the operator 
\[
^{t}M_{f}:\mathcal{S}_{n}^{\prime }\rightarrow \mathcal{S}_{n}^{\prime } 
\]
defined by 
\begin{eqnarray*}
^{t}M_{f}(u)(g) &=&u(M_{f}(g))= \\
&=&u(fg)= \\
&=&fu(g),
\end{eqnarray*}
for every $u$ in $\mathcal{S}_{n}^{\prime }$ and for every $g$ in $\mathcal{S%
}_{n}$. So that, the transpose of the multiplication $M_{f}$ is the
multiplication on $\mathcal{S}_{n}^{\prime }$ by the function $f$. Indeed,
the multiplication of a tempered distribution by an $^{\mathcal{O}_{M}}$%
function is defined by the transpose of $M_{f}$, since this last operator is
self-adjoint with respect to the canonical bilinear form on $\mathcal{S}%
_{n}\times \mathcal{S}_{n}$. In fact, obviously, we have 
\[
\langle M_{f}(g),h\rangle =\langle g,M_{f}(h)\rangle , 
\]
for every pair $(g,h)$ in that Cartesian product $\mathcal{S}_{n}\times 
\mathcal{S}_{n}$. So we can use the standard procedure to extend regular
operators (operators admitting an adjoint with respect to the standard
bilinear form) from their domain $\mathcal{S}_{n}$ to the entire space $%
\mathcal{S}_{n}^{\prime }$.

\bigskip

\subsection{$^{\mathcal{S}}$\textbf{Family of the multiplication operator }$%
M_{f}$}

\bigskip

Since the multiplication operator $M_{f}:\mathcal{S}_{n}\rightarrow \mathcal{%
S}_{n}$ is continuous, we can associate with it an $^{\mathcal{S}}$family $v$%
, in the canonical way. We have 
\begin{eqnarray*}
v_{p} &=&(M_{f}^{\vee })_{p}= \\
&=&\delta _{p}\circ M_{f}= \\
&=&\;^{t}M_{f}(\delta _{p})= \\
&=&f\delta _{p}= \\
&=&f(p)\delta _{p},
\end{eqnarray*}
for every $p$ in $\Bbb{R}^{n}$. In the language of Schwartz matrices we can
say that to the operator $M_{f}$ is associated the Schwartz diagonal matrix $%
f\delta $.

\bigskip

\section{\textbf{Product in }$\mathcal{L}(\mathcal{S}_{n},\mathcal{S}_{m})$%
\textbf{\ by }$^{\mathcal{O}_{M}}$\textbf{functions}}

\bigskip

The basic remark is the following.

\bigskip

\textbf{Proposition.}\emph{\ }$\emph{Let}$\emph{\ }$A\in \mathcal{L}(%
\mathcal{S}_{n},\mathcal{S}_{m})$\emph{\ be a continuous linear operator and
let }$f$\emph{\ be a function of class }$\mathcal{O}_{M}^{(m)}$\emph{. Then,
the mapping} 
\[
fA:\mathcal{S}_{n}\rightarrow \mathcal{S}_{m}:\phi \mapsto fA(\phi ) 
\]
\emph{is a linear and continuous operator too; it is indeed the composition} 
\[
M_{f}\circ A, 
\]
\emph{where }$M_{f}$\emph{\ is the multiplication operator on }$\mathcal{S}%
_{m}$\emph{\ by the function }$f$\emph{.}

\emph{\bigskip }

\emph{Proof.} It is absolutely straightforward. First of all we note that
the product $fA$ is well defined. In fact, we have 
\[
(fA)(\phi )=fA(\phi ), 
\]
and the right-hand function lies in the space $\mathcal{S}_{m}$ because the
function $f$ lies in the space $\mathcal{O}_{M}^{(m\Bbb{)}}$ and the
function $A(\phi )$ lies in the space $\mathcal{S}_{m}$. Moreover, the
bilinear application 
\[
\Phi :\mathcal{O}_{M}^{(m)}\times \mathcal{S}_{m}\rightarrow \mathcal{S}%
_{m}:(f,\psi )\mapsto f\psi 
\]
is separately continuous and we have 
\begin{eqnarray*}
(fA)(\phi ) &=&fA(\phi )= \\
&=&\Phi (f,A(\phi ))= \\
&=&M_{f}(A(\phi )),
\end{eqnarray*}
i.e., 
\begin{eqnarray*}
fA &=&\Phi (f,\cdot )\circ A= \\
&=&M_{f}\circ A,
\end{eqnarray*}
hence the operator $fA$ is the composition of two linear continuous maps and
then it is a linear and continuous operator. $\blacksquare $

\bigskip

\textbf{Definition.}\emph{\ Let }$A\in \mathcal{L}(\mathcal{S}_{n},\mathcal{S%
}_{m})$\emph{\ and }$f\in \mathcal{O}_{M}^{(m\Bbb{)}}$\emph{. The operator} 
\[
fA:\mathcal{S}_{n}\rightarrow \mathcal{S}_{m}:\phi \mapsto fA(\phi ) 
\]
\emph{is called \textbf{the product of the operator} }$A$\emph{\ \textbf{by
the function} }$f$\emph{.}

\bigskip

\textbf{Proposition.}\emph{\ Let }$A,B\in \mathcal{L}(\mathcal{S}_{n},%
\mathcal{S}_{m})$\emph{\ be two continuous linear operators and }$f,g$\emph{%
\ be two functions in }$\mathcal{O}_{M}^{(m\Bbb{)}}.$\emph{\ Then, we have}

\begin{itemize}
\item[\emph{1)}]  \emph{\ }$(f+g)A=fA+gA$\emph{; }$f(A+B)=fA+fB$\emph{; }$1_{%
\Bbb{R}^{m}}A=A$\emph{, where the function }$1_{\Bbb{R}^{m}}$\emph{\ is the
constant function of }$\Bbb{R}^{m}$\emph{\ into }$\Bbb{K}$\emph{\ with value 
}$1$\emph{;}

\item[\emph{2)}]  \emph{\ the map } 
\[
\Phi :\mathcal{O}_{M}^{(m\Bbb{)}}\Bbb{\times }\mathcal{L}(\mathcal{S}_{n},%
\mathcal{S}_{m})\rightarrow \mathcal{L}(\mathcal{S}_{n},\mathcal{S}%
_{m}):(f,A)\mapsto fA 
\]
\emph{is a bilinear map.}
\end{itemize}

\emph{\bigskip }

\emph{Proof. }It's a straightforward computation.\ $\blacksquare $

\bigskip

The above bilinear application is called \emph{multiplication of operators by%
} $\mathcal{O}_{M}$\ \emph{functions}.

\bigskip

\subsection{\textbf{The algebra }$\mathcal{O}_{M}^{(m)}$}

\bigskip

It's easy to see that the algebraic structure $(\mathcal{O}_{M}^{(m\Bbb{)}%
},+,\cdot )$ is a commutative ring with identity, with respect to the usual
pointwise addition and multiplications. For instance, the multiplication is
the operation 
\[
\cdot \;:\mathcal{O}_{M}^{(m\Bbb{)}}\times \mathcal{O}_{M}^{(m\Bbb{)}}\Bbb{\
\rightarrow }\mathcal{O}_{M}^{(m\Bbb{)}}:\Bbb{(}f,g)\mapsto fg, 
\]
where, obviously, if $f,g\in \mathcal{O}_{M}^{(m\Bbb{)}}$, then the
pointwise product $fg$ still lies in $\mathcal{O}_{M}^{(m\Bbb{)}}$. The
identity of the ring is the function $1_{m}:=1_{\Bbb{R}^{m}}$. Moreover, we
have that the subspace $\mathcal{S}_{m}$ of the space $\mathcal{O}_{M}^{(m)}$
is an ideal of the ring $\mathcal{O}_{M}^{(m)}$. The subring of $\mathcal{O}%
_{M}^{(m)}$ formed by the invertible elements of $\mathcal{O}_{M}^{(m)}$ is
exactly the multiplicative subgroup of those elements $f$ such that the
multiplicative inverse $f^{-1}$ belongs to the space $\mathcal{O}_{M}^{(m)}$
too.

\bigskip

So that, the space $\mathcal{O}_{M}^{(m)}$ is a locally convex topological
algebra with unit element.

\bigskip

\subsection{\textbf{The module }$\mathcal{L}(\mathcal{S}_{n},\mathcal{S}
_{m}) $}

\bigskip

\textbf{Proposition.}\emph{\ Let }$\cdot $\emph{\ be the multiplication by }$%
\mathcal{O}_{M}^{(m)}$\emph{\ functions\ defined in the above theorem. Then,
the algebraic structure }$(\mathcal{L}(\mathcal{S}_{n},\mathcal{S}%
_{m}),+,\cdot )$\emph{\ is a left module over the ring }$(\mathcal{O}_{M}^{(m%
\Bbb{)}},+,\cdot ).$

\emph{\bigskip }

\emph{Proof. }Recalling the preceding theorem, we have to prove only the
pseudo-associative law, i.e. we have to prove that for every couple of
functions $f,g\in \mathcal{O}_{M}^{(m\Bbb{)}}$ and for every linear
continuous operator $A\in \mathcal{L}(\mathcal{S}_{n},\mathcal{S}_{m})$, we
have 
\[
(fg)A=f(gA). 
\]
In fact, for each $\phi \in \mathcal{S}_{n}$, we have 
\begin{eqnarray*}
\lbrack (fg)A](\phi ) &=&(fg)A(\phi )= \\
&=&f(gA(\phi ))= \\
&=&f(gA)(\phi ))= \\
&=&[f(gA)](\phi ),
\end{eqnarray*}
as we desired. $\blacksquare $

\bigskip

\section{\textbf{Products of }$^{\mathcal{S}}$\textbf{families by }$^{%
\mathcal{O}_{M}}$\textbf{functions}}

\bigskip

The central definition of the chapter is the following.

\bigskip

\textbf{Definition (product of Schwartz families by smooth functions). }%
\emph{Let }$v\in \mathcal{S}(\Bbb{R}^{m},\mathcal{S}_{n}^{\prime })$\emph{\
be an }$^{\mathcal{S}}$\emph{family of distributions and let }$f\in \mathcal{%
C}^{\infty }(\Bbb{R}^{m},\Bbb{K)}$\emph{\ be a smooth function. The \textbf{%
product} \textbf{of the family }}$v$\emph{\ \textbf{by the function} }$f$%
\textbf{\ }\emph{is the family} 
\[
fv:=(f(p)v_{p})_{p\in \Bbb{R}^{m}}. 
\]

\bigskip

\textbf{Theorem.}\emph{\ Let }$v\in \mathcal{S}(\Bbb{R}^{m},\mathcal{S}%
_{n}^{\prime })$\emph{\ be an }$^{\mathcal{S}}$\emph{family and }$f\in 
\mathcal{O}_{M}^{(m\Bbb{)}}$\emph{. Then, the family }$fv$\emph{\ lies in }$%
\mathcal{S}(\Bbb{R}^{m},\mathcal{S}_{n}^{\prime }).$\emph{\ Moreover, we have%
} 
\[
(fv)^{\wedge }=f\widehat{v}. 
\]
\emph{Consequently, concerning the superposition operator of the family }$fv$%
\emph{, since }$f\widehat{v}=M_{f}\circ \widehat{v}$\emph{, we have} 
\[
^{t}(fv)^{\wedge }=\;^{t}\widehat{v}\circ \;^{t}M_{f}, 
\]
\emph{or equivalently, in superposition form} 
\[
\int_{\Bbb{R}^{m}}a(fv)=\int_{\Bbb{R}^{m}}(fa)v, 
\]
\emph{for every coefficient distribution }$a$\emph{\ in }$\mathcal{S}%
_{m}^{\prime }$\emph{.}

\emph{\bigskip }

\emph{Proof.} Let $\phi \in \mathcal{S}_{n}$ be a test function, we have 
\begin{eqnarray*}
(fv)(\phi )(p) &=&(fv)_{p}(\phi )= \\
&=&(f(p)v_{p})(\phi )= \\
&=&f(p)v_{p}(\phi )= \\
&=&f(p)\widehat{v}(\phi )(p)
\end{eqnarray*}
and hence the function $(fv)(\phi )$ equals $f\widehat{v}(\phi )$, which
lies in $\mathcal{S}_{m}$. Thus, the product $fv$ lies in the space of
Schwartz families $\mathcal{S}(\Bbb{R}^{m},\mathcal{S}_{n}^{\prime })$. For
any test function $\phi \in \mathcal{S}_{n}$, by the above consideration, we
deduce 
\[
(fv)^{\wedge }(\phi )=f\widehat{v}(\phi ), 
\]
that is, the equality of operators 
\[
(fv)^{\wedge }=f\widehat{v}, 
\]
where $f\widehat{v}$ is the product of the operator $\widehat{v}$ by the
function $f$, product which belongs to the space $\mathcal{L}(\mathcal{S}%
_{m},\mathcal{S}_{n})$. Moreover, concerning the superposition operator of
the family $fv$, we obtain 
\begin{eqnarray*}
\int_{\Bbb{R}^{m}}a(fv) &=&\;^{t}(fv)^{\wedge }(a)= \\
&=&\;^{t}(f\widehat{v})(a)= \\
&=&\;^{t}(M_{f}\circ \widehat{v})(a)= \\
&=&\;(^{t}\widehat{v}\circ \;^{t}M_{f})(a)= \\
&=&\;^{t}\widehat{v}(^{t}M_{f}(a))= \\
&=&\;^{t}\widehat{v}(fa)= \\
&=&\int_{\Bbb{R}^{m}}(fa)v,
\end{eqnarray*}
for every distribution $a$ in $\mathcal{S}_{m}^{\prime }$. $\blacksquare $

\bigskip

\textbf{Theorem.}\emph{\ Let }$f,g$\emph{\ two functions in the space }$%
\mathcal{O}_{M}^{(m\Bbb{)}}$\emph{\ and }$v,w$\emph{\ two Schwartz families
in the space }$\mathcal{S}(\Bbb{R}^{m},\mathcal{S}_{n}^{\prime }).$\emph{\
Then, we have:}

\begin{itemize}
\item[\emph{1)}]  \emph{\ }$(f+g)v=fv+gv,\;f(v+w)=fv+fw$\emph{\ and }$%
1_{m}v=v$\emph{;}

\item[\emph{2)}]  \emph{\ the map } 
\[
\Phi :\mathcal{O}_{M}^{(m\Bbb{)}}\times \mathcal{S}(\Bbb{R}^{m},\mathcal{S}%
_{n}^{\prime })\rightarrow \mathcal{S}(\Bbb{R}^{m},\mathcal{S}_{n}^{\prime
}):(f,v)\mapsto fv 
\]
\emph{is a bilinear map.}
\end{itemize}

\emph{\bigskip }

\emph{Proof. }1) For all $p\in \Bbb{R}^{m}$, we have 
\begin{eqnarray*}
\left[ \left( f+g\right) v\right] (p) &=&(f+g)(p)v_{p}= \\
&=&(f(p)+g(p))v_{p}= \\
&=&f(p)v_{p}+g(p)v_{p}= \\
&=&(fv)_{p}+(gv)_{p},
\end{eqnarray*}
i.e. $(f+g)v=fv+gv$. For all $p\in \Bbb{R}^{m}$, we have 
\begin{eqnarray*}
\left[ f\left( v+w\right) \right] (p) &=&f(p)(v+w)_{p}= \\
&=&f(p)(v_{p}+w_{p})= \\
&=&f(p)v_{p}+f(p)w_{p}= \\
&=&(fv)_{p}+(fw)_{p},
\end{eqnarray*}
i.e. $f(v+w)=fv+fw$. For all $p\in \Bbb{R}^{m}$, we have 
\[
(1_{\Bbb{R}^{m}}v)(p)=1_{\Bbb{R}^{m}}(p)v_{p}=v_{p}; 
\]
i.e. $1_{\Bbb{R}^{m}}v=v$. 2) follows immediately by 1). $\blacksquare $

\bigskip

The bilinear application of the point 2) of the preceding theorem is called 
\emph{multiplication of Schwartz families by }$\mathcal{O}_{M}$\emph{\
functions.}

\bigskip

\textbf{Theorem (of structure). }\emph{Let }$\cdot $\emph{\ the operation} 
\emph{defined above. Then, the algebraic structure} $(\mathcal{S}(\Bbb{R}
^{m},\mathcal{S}_{n}^{\prime }),+,\cdot )$\emph{\ is a left module over the
ring }$(\mathcal{O}_{M}^{(m\Bbb{)}},+,\cdot ).$

\bigskip

\emph{Proof. }It's analogous to the proof of the corresponding proposition
for operators. $\blacksquare $

\bigskip

\textbf{Theorem (of isomorphism). }\emph{The application} 
\[
(\cdot )^{\wedge }:\mathcal{S}(\Bbb{R}^{m},\mathcal{S}_{n}^{\prime
})\rightarrow \mathcal{L}(\mathcal{S}_{n},\mathcal{S}_{m}) 
\]
\emph{is a module isomorphism.}

$\bigskip $

$\emph{Proof.}$ It follows easily from the above theorem. $\blacksquare $

\bigskip

\section{$^{\mathcal{O}_{M}}$\textbf{Functions and Schwartz basis}}

\bigskip

In this section we study some important relations among a Schwartz family $w$
and its multiples $fw$.

\bigskip

\textbf{Theorem.}\emph{\ Let }$w\in \mathcal{S}(\Bbb{R}^{m},\mathcal{S}%
_{n}^{\prime })$\emph{\ be a Schwartz family and let }$f\in \mathcal{O}%
_{M}^{(m\Bbb{)}}$\emph{. Then, the hull }$^{\mathcal{S}}\mathrm{span}(w)$%
\emph{\ of the family }$w$ \emph{contains the hull }$^{\mathcal{S}}\mathrm{%
span}(fw)$\emph{\ of the multiple family }$fw$\emph{. Moreover, if a
distribution }$a\;$\emph{represents the distribution }$u$\emph{\ with
respect to the family }$fw$\emph{\ (that is, if }$u=a.(fw)$\emph{) then the
distribution }$fa$\emph{\ represents the distribution }$u$\emph{\ with
respect to the family }$w$ \emph{(that is, if }$u=(fa).w$\emph{).}

\emph{\bigskip }

\emph{Proof.} 1) Let $u$ be a vector of the $^{\mathcal{S}}$linear hull $^{%
\mathcal{S}}\mathrm{span}(fw)$. Then, there exists a coefficient
distribution $a\in \mathcal{S}_{m}^{\prime }$ such that 
\[
u=\int_{\Bbb{R}^{m}}a(fw), 
\]
and this is equivalent (as we already have seen) to the equality 
\[
u=\int_{\Bbb{R}^{m}}(fa)w; 
\]
hence the vector $u$ belongs also to the $^{\mathcal{S}}$linear hull $^{%
\mathcal{S}}\mathrm{span}(w)$. Hence the $^{\mathcal{S}}$linear hull $^{%
\mathcal{S}}\mathrm{span}(fw)$ is contained in the $^{\mathcal{S}}$linear
hull $^{\mathcal{S}}\mathrm{span}(w)$. $\blacksquare $

\bigskip

\textbf{Theorem.}\emph{\ Let }$w\in \mathcal{S}(\Bbb{R}^{m},\mathcal{S}%
_{n}^{\prime })$\emph{\ be a Schwartz family and let }$f\in \mathcal{O}%
_{M}^{(m\Bbb{)}}$\emph{\ be a function different from }$0$\emph{\ at every
point of its domain. Then, the following assertions hold true:}

\begin{itemize}
\item[\emph{1)}]  \emph{\ if the family }$w$\emph{\ is }$^{\mathcal{S}}$%
\emph{linearly independent, the family }$fw$\emph{\ is }$^{\mathcal{S}}$%
\emph{linearly independent too;}

\item[\emph{2)}]  \emph{\ the Schwartz linear hull }$^{\mathcal{S}}\mathrm{%
span}(w)$\emph{\ contains the hull }$^{\mathcal{S}}\mathrm{span}(fw)$\emph{;}

\item[\emph{3)}]  \emph{\ if the family }$w\;$\emph{is }$^{\mathcal{S}}$%
\emph{linearly independent, for each vector }$u$\emph{\ in the hull }$^{%
\mathcal{S}}\mathrm{span}(fw)$\emph{, we have } 
\[
\lbrack u\mid w]=f[u\mid fw], 
\]
\emph{where, as usual, by }$[u|v]$\emph{\ we denote the Schwartz coordinate
system of a distribution }$u$\emph{\ (in the Schwartz linear hull of }$v$%
\emph{) with respect to a Schwartz linear independent family }$v$\emph{;}

\item[\emph{4)}]  \emph{\ if the family }$w$\emph{\ is an }$^{\mathcal{S}}$%
\emph{basis of a subspace }$V$\emph{, then }$fw$\emph{\ is an }$^{\mathcal{S}%
}$\emph{basis of its }$^{\mathcal{S}}$\emph{linear hull }$^{\mathcal{S}}%
\mathrm{span}(fw)$\emph{\ (that in general is a proper subspace of the hull }%
$^{\mathcal{S}}\mathrm{span}(w)$\emph{).}
\end{itemize}

\emph{\bigskip }

\emph{Proof.} 1) Let $a\in \mathcal{S}_{m}^{\prime }$ be such that 
\[
\int_{\Bbb{R}^{m}}a(fw)=0_{\mathcal{S}_{n}^{\prime }}, 
\]
we have 
\begin{eqnarray*}
0_{\mathcal{S}_{n}^{\prime }} &=&\int_{\Bbb{R}^{m}}a(fw)= \\
&=&\int_{\Bbb{R}^{m}}(fa)w,
\end{eqnarray*}
thus, because the family $w$ is $^{\mathcal{S}}$linearly independent we have 
$fa=0_{\mathcal{S}_{n}^{\prime }}$. Since $f$ is different from $0$ at every
point, we can conclude $a=0_{\mathcal{S}_{n}^{\prime }}$.

2) Let $u$ be a vector of the Schwartz linear hull $^{\mathcal{S}}\mathrm{%
span}(fw)$. Then, there exists a coefficient distribution $a\in \mathcal{S}%
_{m}^{\prime }$ such that 
\[
u=\int_{\Bbb{R}^{m}}a(fw), 
\]
or equivalently such that 
\[
u=\int_{\Bbb{R}^{m}}(fa)w, 
\]
and hence the vector $u$ belongs also to the hull $^{\mathcal{S}}\mathrm{span%
}(w)$. Hence the Schwartz linear hull $^{\mathcal{S}}\mathrm{span}(fw)$ is
contained in the hull $^{\mathcal{S}}\mathrm{span}(w)$.

3) If the family $w$ is $^{\mathcal{S}}$linearly independent, from the above
two equalities, we deduce $(u)_{fw}=a$ and $(u)_{w}=fa$, from which 
\begin{eqnarray*}
(u)_{w} &=&fa= \\
&=&f(u)_{fw},
\end{eqnarray*}
as we claimed.

4) is an obvious consequence of the preceding properties. $\blacksquare $

\bigskip

\section{$^{\mathcal{O}_{M}}$\textbf{Invertible functions and }$^{\mathcal{S}
}$\textbf{basis}}

\bigskip

We recall that an invertible element of $\mathcal{O}_{M}^{(m)}$ is any
function $f$ everywhere different from $0$ and such that its multiplicative
inverse $f^{-1}$ lives in $\mathcal{O}_{M}^{(m)}$ too. The set of the
invertible elements of the space $\mathcal{O}_{M}^{(m)}$ is a group with
respect to the pointwise multiplication, and we will denote it by $\mathcal{G%
}_{M}^{(m)}$.

\bigskip

\textbf{Theorem.}\emph{\ Let }$w\in \mathcal{S}(\Bbb{R}^{m},\mathcal{S}%
_{n}^{\prime })$\emph{\ be a Schwartz family and let }$f\in \mathcal{G}%
_{M}^{(m)}$\emph{\ be an invertible element of the ring }$\mathcal{O}%
_{M}^{(m)}$\emph{\ (in particular, it must be a function different form }$0$%
\emph{\ at every point). Then, the following assertions hold true:}

\begin{itemize}
\item[\emph{1)}]  \emph{\ the family }$w$\emph{\ is }$^{\mathcal{S}}$\emph{%
linearly independent if and only if the multiple family }$fw$\emph{\ is }$^{%
\mathcal{S}}$\emph{linearly independent;}

\item[\emph{2)}]  \emph{\ the hull }$^{\mathcal{S}}\mathrm{span}(w)$\emph{\
coincides with the hull }$^{\mathcal{S}}\mathrm{span}(fw)$\emph{;}

\item[\emph{3)}]  \emph{\ if the family }$w\;$\emph{is }$^{\mathcal{S}}$%
\emph{linearly independent, then, for each vector }$u$\emph{\ in the hull }$%
^{\mathcal{S}}\mathrm{span}(w)$\emph{, we have} 
\[
\lbrack u\mid fw]=\left( 1/f\right) [u\mid w], 
\]
\emph{where, as usual, by }$[u|v]$\emph{\ we denote the Schwartz coordinate
system of a distribution }$u$\emph{\ (in the Schwartz linear hull of }$v$%
\emph{) with respect to a Schwartz linear independent family }$v$\emph{;}

\item[\emph{4)}]  \emph{\ the family }$w$\emph{\ is an }$^{\mathcal{S}}$%
\emph{basis of a subspace }$V$\emph{\ if and only if its multiple }$fw$\emph{%
\ is an }$^{\mathcal{S}}$\emph{basis of the }$^{\mathcal{S}}$\emph{linear
hull }$^{\mathcal{S}}\mathrm{span}(fw)$\emph{\ (that in this case coincides
with }$^{\mathcal{S}}\mathrm{span}(w)$\emph{).}
\end{itemize}

\bigskip

\emph{Proof.} 1) Let $a\in \mathcal{S}_{m}^{\prime }$ be a distribution such
that 
\[
\int_{\Bbb{R}^{m}}aw=0_{\mathcal{S}_{n}^{\prime }}, 
\]
we have 
\begin{eqnarray*}
0_{\mathcal{S}_{n}^{\prime }} &=&\int_{\Bbb{R}^{m}}aw= \\
&=&\int_{\Bbb{R}^{m}}(f^{-1}a)(fw),
\end{eqnarray*}
thus, because $fw$ is $^{\mathcal{S}}$linearly independent we have $%
f^{-1}a=0_{\mathcal{S}_{n}^{\prime }}$. Since $f^{-1}$ is different form $0$
at every point we can conclude $a=0_{\mathcal{S}_{n}^{\prime }}$.

2) Let $u$ be in $^{\mathcal{S}}\mathrm{span}(w)$. Then, there exists a
distribution $a\in \mathcal{S}_{m}^{\prime }$ such that 
\[
u=\int_{\Bbb{R}^{m}}aw. 
\]
Now, we have 
\[
u=\int_{\Bbb{R}^{m}}(f^{-1}a)\left( fw\right) , 
\]
so the distribution $u$ lies in $^{\mathcal{S}}\mathrm{span}(fw)$, and hence 
$^{\mathcal{S}}\mathrm{span}(w)$ is contained in $^{\mathcal{S}}\mathrm{span}%
(fw)$. Vice versa, let $u$ be in $^{\mathcal{S}}\mathrm{span}(fw)$. Then,
there exists a distribution $a\in \mathcal{S}_{m}^{\prime }$ such that 
\[
u=\int_{\Bbb{R}^{m}}a(fw). 
\]
Now, we have (equivalently) 
\[
u=\int_{\Bbb{R}^{m}}(fa)w, 
\]
and hence $u$ lies also in $^{\mathcal{S}}\mathrm{span}(w)$, hence $^{%
\mathcal{S}}\mathrm{span}(fw)$ is contained in $^{\mathcal{S}}\mathrm{span}%
(w)$ (as we already have seen in the general case). Concluding 
\[
^{\mathcal{S}}\mathrm{span}(w)=\;^{\mathcal{S}}\mathrm{span}(fw). 
\]

3) For any distribution $u$ in the Schwartz linear hull of the family $w$,
we have 
\[
u=\int_{\Bbb{R}^{m}}\left[ u\mid w\right] w, 
\]
hence 
\[
u=\int_{\Bbb{R}^{m}}(f^{-1}[u|w])\left( fw\right) , 
\]
as we desired.

4) It follows immediately from the above properties. $\blacksquare $

\bigskip

\textbf{Theorem.}\emph{\ Let }$e\in \mathcal{B}(\Bbb{R}^{m},\mathcal{S}%
_{n}^{\prime })$\emph{\ be an }$^{\mathcal{S}}$\emph{basis of the space }$%
\mathcal{S}_{n}^{\prime }$\emph{\ and let }$f\in \mathcal{O}_{M}^{(m)}$\emph{%
. Then the multiple }$fe$\emph{\ is an }$^{\mathcal{S}}$\emph{basis of the
space }$\mathcal{S}_{n}^{\prime }$ \emph{if and only if the factor }$f$\emph{%
\ is an invertible element of the ring }$\mathcal{O}_{M}^{(m)}$\emph{.}

\emph{\bigskip }

\emph{Proof.} We must prove that, if $fe$ is an $^{\mathcal{S}}$basis of $%
\mathcal{S}_{n}^{\prime }$, then $f$ is an invertible element of the ring $%
\mathcal{O}_{M}^{(m)}$. First of all observe that, since $fe$ is a basis,
then $fe$ is $^{\mathcal{S}}$linearly independent and consequently linearly
independent in the ordinary algebraic sense; consequently every distribution 
$f(p)e_{p}$ must be a non zero distribution and this implies that any value $%
f(p)$ must be different from $0$, so we can consider the multiplicative
inverse $f^{-1}$. We now have to prove that the multiplicative inverse $%
f^{-1}$ lives in $\mathcal{O}_{M}^{(m)}$, or equivalently that, for every
test function $g$ in $\mathcal{S}_{m}$, the product $f^{-1}g$ lives in $%
\mathcal{S}_{m}$. For, let $g$ be in $\mathcal{S}_{m}$, since $fe$ is a
basis, its associated operator from $\mathcal{S}_{n}$ into $\mathcal{S}_{m}$
is surjective, then there is a function $h$ in $\mathcal{S}_{n}$ such that $%
(fe)^{\wedge }(h)=g$, the last equality is equivalent to 
\[
fe(h)=g, 
\]
that is 
\[
f^{-1}g=e(h), 
\]
so that $f^{-1}g$ actually lives in the space $\mathcal{S}_{m}$. $%
\blacksquare $

\bigskip

We can generalize the above result as it follows.

\bigskip

\textbf{Theorem.}\emph{\ Let }$e\in \mathcal{B}(\Bbb{R}^{m},V)$\emph{\ be an 
}$^{\mathcal{S}}$\emph{basis of a (weakly*) closed subspace }$V$\emph{\ of
the space }$\mathcal{S}_{n}^{\prime }$\emph{\ and let }$f\in \mathcal{O}
_{M}^{(m)}$\emph{. Then the multiple family }$fe$\emph{\ is an }$^{\mathcal{S%
}}$\emph{basis of the subspace }$V$ \emph{if and only if the factor }$f$%
\emph{\ is an invertible element of the ring }$\mathcal{O}_{M}^{(m)}$\emph{.}

\emph{\bigskip }

\emph{Proof.} We must prove that, if $fe$ is an $^{\mathcal{S}}$basis of the
subspace $V$, then $f$ is an invertible element of the ring $\mathcal{O}%
_{M}^{(m)}$. First of all observe that, since $fe$ is a basis, then $fe$ is $%
^{\mathcal{S}}$linearly independent and consequently linearly independent in
the ordinary algebraic sense; consequently every distribution $f(p)e_{p}$
must be a non zero distribution and this implies that any value $f(p)$ must
be different from $0$. So we can consider its multiplicative inverse $f^{-1}$%
. We now have to prove that the multiplicative inverse $f^{-1}$ lives in the
space $\mathcal{O}_{M}^{(m)}$, or equivalently that, for every test function 
$g$ in $\mathcal{S}_{m}$ the product $f^{-1}g$ lives in $\mathcal{S}_{m}$.
For, let $g$ be in $\mathcal{S}_{m}$, since $fe$ is an $^{\mathcal{S}}$basis
of the topologically closed subspace $V$, its associated operator $%
(fe)^{\wedge }$ from $\mathcal{S}_{n}$ into $\mathcal{S}_{m}$ is surjective
(this follows, by the closedness of $V$, from the Dieudonn\'{e}-Schwartz
theorem, since the transpose of the operator $(fe)^{\wedge }$ is the
superposition operator of $fe$, which is injective since the family $fe$ is
Schwartz linearly independent). Hence, by surjectivity, there is a function $%
h$ in $\mathcal{S}_{n}$ such that $(fe)^{\wedge }(h)=g$, the last equality
is equivalent to the following one 
\[
fe(h)=g, 
\]
that is 
\[
f^{-1}g=e(h), 
\]
so that the function $f^{-1}g$ actually lives in the space $\mathcal{S}_{m}$%
. $\blacksquare $

\bigskip

\textbf{David Carf\`{i}}

\emph{Faculty of Economics}

\emph{University of Messina}

\emph{davidcarfi71@yahoo.it}

\end{document}